\documentclass[a4paper,11pt]{article}
\usepackage{amsthm}

\usepackage[T2A]{fontenc}
\usepackage[cp1251]{inputenc}
\usepackage[english,russian]{babel}
\usepackage[tbtags]{amsmath}
\usepackage{amsfonts,amssymb}

\usepackage{mathrsfs}
\usepackage{graphicx}

\usepackage{euscript,textcomp,verbatim,fancyhdr}
\usepackage[all]{xy}
\usepackage{epsfig}
\usepackage{latexsym}
\usepackage{mdwtab}           
\usepackage{graphics}

\theoremstyle{plain}

\newtheorem{Thm}{Теорема}[subsection]
\newtheorem{Lem}[Thm]{Лемма}
\newtheorem{Sta}[Thm]{Утверждение}
\newtheorem{Cor}[Thm]{Следствие}
\theoremstyle{definition}

\newtheorem{Def}[Thm]{Определение}
\newtheorem{Not}[Thm]{Обозначение}

\newcommand{\const}{\mathrm{const}}

\newcommand{\lam}{\lambda}
\newcommand{\RR}{\mathbb R}

\newcommand{\ZZ}{\mathbb Z}

\newcommand{\FF}{\mathbb{F}}
\newcommand{\XX}{\mathbb{\widetilde X}}
\newcommand{\YY}{\mathbb{\widetilde Y}}
\newcommand{\uups}{\boldsymbol{\ups}}
\newcommand{\barc}{{\boldsymbol{c}}}

\newcommand{\QQ}{\mathbb{Q}}
\newcommand{\UU}{\mathbb{U}}
\newcommand{\KK}{\mathbb{\widetilde K}}

\newcommand{\T}{{\cal T}}
\renewcommand{\tilde}{\widetilde}
\renewcommand{\hat}{\widehat}

\newcommand{\MM}{{\widetilde{\cal M}}}
\newcommand{\IMM}{{\cal M}}
\newcommand{\D}{{\mathscr D}}
\newcommand{\N}{{\cal C}}
\newcommand{\Diff}{{\rm Diff}}

\newcommand{\Hom}{{\rm Hom}}

\renewcommand{\Im}{{\rm Im}}
\newcommand{\num}{{\rm num}}

\newcommand{\diskr}{{\rm discr}}
\newcommand{\Forg}{{\rm Forg}}
\newcommand{\Ev}{{\rm Ev}}

\newcommand{\eps}{\varepsilon}
\newcommand{\ups}{\upsilon}
\newcommand{\Int}{{\rm int\,}}
\newcommand{\grad}{{\rm grad\,}}
\newcommand{\stab}{{\rm stab}}
\newcommand{\st}{{\rm st}}

\newcommand{\id}{{\rm id}}
\newcommand{\isot}{{\rm isot}}
\newcommand{\glue}{{\rm glue}}

\renewcommand{\mod}{{\rm \,mod\,}}

\renewcommand{\:}{\colon\,}
\renewcommand{\a}{x}
\renewcommand{\b}{y}

\renewcommand{\d}{\partial}

\renewcommand{\i}{{\tilde s}}

\newcommand{\aapprox}{\sim} 
\renewcommand{\emptyset}{\varnothing}

\begin{document}


\date{}

\author{E.\,A.~Kudryavtseva}

\title{Topology of the spaces of Morse functions on surfaces
}


\maketitle


Let $M$ be a smooth closed orientable surface, and let $F$ be the space
of Morse functions on $M$ such that at least $\chi(M)+1$ critical points of each function of $F$ are labeled
by different labels (enumerated). Endow the space $F$ with $C^\infty$-topology.
We prove the homotopy equivalence $F\sim R\times\MM$ where $R$ is one of the manifolds $\RR P^3$, $S^1\times S^1$ and the point in dependence on the sign of $\chi(M)$, and $\MM$ is the universal moduli space of framed Morse functions, which is a smooth stratified manifold. Morse inequalities for the Betti numbers of the space $F$ are obtained.


{\bf MSC-class:} 58E05, 57M50, 58K65, 46M18

\bigskip
УДК 515.164.174+515.164.22+515.122.55

\begin{center}
{\LARGE Топология пространств функций Морса на поверхностях} \bigskip\\
\large Е.\,А.~Кудрявцева
\bigskip
\end{center}

\begin{abstract}
Пусть $M$ --- гладкая замкнутая ориентируемая поверхность, и пусть $F$
-- пространство функций Морса на $M$, у которых не менее чем $\chi(M)+1$
критических точек помечены различными метками (пронумерованы).
Снабдим $C^\infty$-топологией пространство $F$.
Доказана гомотопическая эквивалентность $F\sim R\times\MM$, где $R$ --- одно
из многообразий $\RR P^3$, $S^1\times S^1$ и точка в зависимости от
знака $\chi(M)$, а $\MM$ -- универсальное пространство модулей оснащенных функций Морса, являющееся гладким стратифицированным многообразием. Получены неравенства Морса для чисел Бетти пространства $F$.


\medskip
{\bf Ключевые слова:} функция Морса, оснащенная функция Морса,
комплекс оснащенных функций Морса, $C^\infty$-топология,
универсальное пространство модулей.
\end{abstract}

\footnotetext{Работа выполнена при поддержке РФФИ (грант
\No~10–01–00748-а), Программы поддержки ведущих научных школ РФ
(грант \No~НШ-3224.2010.1), Программы ``Развитие научного потенциала
высшей школы'' (грант \No~2.1.1.3704), ФЦП ``Научные и
научно-педагогические кадры инновационной России'' (гранты
\No~02.740.11.5213 и \No~14.740.11.0794).}


\subsection{Введение}\label{sec:intro}

В настоящей работе изучается топология пространства $F=F(M)$ функций
Морса на компактной гладкой двумерной поверхности $M$.
Предполагается, что у каждой функции $f\in F$ по меньшей мере
$\chi(M)+1$ критических точек помечены различными метками
(пронумерованы). В работе \cite {kp1} введено понятие оснащенной
функции Морса (см.\ определение \ref {def:framed}) и доказана
гомотопическая эквивалентность $F\sim\FF$ пространства $F$ функций
Морса и пространства $\FF=\FF(M)$ оснащенных функций Морса
(\cite{kp1,KP2}). В работе \cite {KP3} построены комплекс $\KK$
оснащенных функций Морса и содержащее его гладкое стратифицированное
многообразие $\MM$ (см.\ утвер\-жде\-ние \ref {sta:MM1}). Мы
доказываем (теорема~\ref {thm:KP4add}), что пространство $F$ функций
Морса гомотопически эквивалентно полиэдру $R\times\MM$, где $R=R(M)$
-- одно из многообразий $\RR P^3$, $S^1$, $S^1\times S^1$ и точка
(см.\ (\ref {eq:EE})). Таким образом, наш результат сводит изучение
топологии пространства $F$ функций Морса к комбинаторной задаче ---
изучению топологии многообразия $\MM$. Гомологии многообразия $\MM$
могут быть изучены с помощью его естественной стратификации, а также
индуцированной стратификации специальной окрестности
$\MM_{\succeq[f]_\isot}$ каждого страта $\MM_{[f]_\isot}$ (см.\
утверждение~\ref {sta:MM1}). Этим методом мы получаем в случае
$M=S^2$ неравенства Морса для чисел Бетти многообразия $\MM$ и
находим его эйлерову характеристику (следствие~\ref {cor:ineq}).

Вопросы о линейной связности пространств функций Морса на поверхности изучались
С.В.\ Матвеевым~\cite{Kmsb}, Х.~Цишангом, В.В.\ Шарко~\cite{SH},
Е.А.\ Кудрявцевой~\cite{Kmsb}, С.И.\ Макси\-мен\-ко \cite{Max2005},
Ю.М.\ Бурманом~\cite{Bu,Bu2} (для пространств гладких функций без
критических точек на открытых поверхностях) и Е.А.\
Кудрявцевой~\cite{BaK}
(для пространств функций Морса с фиксированным множеством критических точек).
Количество классов эквивалентности (см.\ определение \ref {def:equiv}) простых функций Морса на поверхности
исследовалось в работе \cite {Kul}, а топология классов изотопности -- в работе
\cite{Max}. В работах \cite {F86,FZ0,BF0,BFbook,K,K:stab,BrK:stab} функции Морса
изучались в связи с задачей классификации (лиувиллевой, орбитальной)
не\-вы\-рожд\-ен\-ных интегри\-ру\-е\-мых гамильтоновых систем с
двумя степенями свободы. Группы гомологий и гомотопий пространств
функций с умеренными особенностями (с допущением не\-мор\-сов\-ских
особенностей) на окружности изучался В.И.\ Арнольдом~\cite{A89}.
Невыполнение 1-параметрического $h$-принципа для пространств функций Морса на некоторых компактных многообразиях размерности большей $5$ показано в работах~\cite {ChL09,Hatcher75} (см.\ также \cite [\S1]{kp1}).

Статья имеет следующую структуру. В~\S \ref {sec:analyt}
формулируются основные результаты настоящей работы (теорема~\ref
{thm:KP4add} и следствие \ref {cor:ineq}).
В \S\ref {sec:KK0} описывается конструкция из \cite {KP3} гладкого
страти\-фи\-цированного $3q$-мерного многообразия $\MM$,
где $q$ -- количество седловых кри\-ти\-чес\-ких точек функций Морса
из $F$ (см.\ определение~\ref {def:Morse} и утверждение~\ref
{sta:MM1}).
В \S\ref {sec:MM} доказывается, что многообразие $\MM$ гомеоморфно
универсальному пространству модулей $\FF^1/\D^0$ оснащенных функций
Морса (утверждение \ref {sta:MM2}). В \S\ref {sec:FF:MM}
уста\-нав\-ли\-ва\-ет\-ся гомеоморфизм $\FF^1\approx\D^0\times\MM$
(утверждение \ref {sta:MM3}).

Автор приносит благодарность С.А.\ Мелихову, Д.А.~Пермякову и А.Т.~Фоменко за внимание к работе и полезные
обсуждения.

\subsection{Основные понятия и формулировка основных результатов} \label {sec:analyt}

 \begin{Def} [обобщенное пространство функций Морса {\rm\cite {kp1}}] \label {def:Morse}
Пусть $M$ --- гладкая (т.е.\ класса $C^\infty$) компактная связная
поверхность, край которой пуст или не пуст, $\d M=\d^+M\sqcup\d^-M$,
где $\d^+M$ --- объединение некоторых граничных окружностей.
Пусть $d^+,d^-\ge0$ --- число окружностей в $\d^+M$ и $\d^-M$ соответственно.

{\rm(A)} Обозначим через $C^\infty(M)$ пространство гладких (т.е.\
класса $C^\infty$) вещест\-вен\-но\-знач\-ных функций $f$ на $M$.
Обозначим через $C^\infty(M,\d^+M,\d^-M)\subset C^\infty(M)$
подпространство, состоящее из таких функций $f\in C^\infty(M)$, что
все ее критические точки (т.е.\ такие точки $x\in M$, что $df|_x=0$)
принадлежат $\Int M$, а любая граничная точка $x\in\d M$ имеет такую
окрестность $U$ в $M$, что $f(U\cap\d M)=f(x)$, причем
$\inf(f|_U)=f(x)$ при $x\in\d^-M$, и $\sup(f|_U)=f(x)$ при
$x\in\d^+M$.

{\rm(B)}
Пусть $\tilde F:=F_{p,q,r}(M,\d^+M,\d^-M)$ --- пространство функций
Морса $f$ на по\-верх\-нос\-ти \\ $(M,\d^+M,\d^-M)$, имеющих
ровно $p$ критических точек локальных ми\-ни\-му\-мов, $q$ седловых
точек и $r$ точек локальных максимумов.
Обозначим через $F^\num$ про\-стран\-ство, полученное из $\tilde F$
введением нумерации у некоторых из критических точек (называемых
отмеченными) для функций Морса $f\in\tilde F$.
Обозначим ко\-ли\-чест\-во отмеченных критических точек локальных
минимумов, максимумов и седловых точек через $\hat p,\hat r,\hat q$
соответственно, $0\le\hat p\le p$, $0\le\hat q\le q$, $0\le\hat r\le
r$.

{\rm(C)} Пусть $0\le p^*\le\hat p$, $0\le q^*\le\hat q$, $0\le r^*\le\hat r$. Обозначим
 $$
(p',p'';q',q'';r',r''):=(\hat p-p^*,p-\hat p;\hat q-q^*,q-\hat q;\hat
r-r^*,r-\hat r).
 $$
Для каждой функции $f\in F^\num$ обозначим через $\N_{f,\lam}$
множество ее критических точек индекса $\lam$,
и через $\hat\N_{f,\lam}\subseteq\N_{f,\lam}$ множество отмеченных
критических точек, $\lam=0,1,2$. В множестве отмеченных (а потому
занумерованных) критических точек обозначим через $\N^*_{f,0}$,
$\N^*_{f,1}$, $\N^*_{f,2}$ подмножество, состоящее из первых
$p^*,q^*,r^*$ точек соответственно. Фиксируем ``базисную'' функцию $f_*\in F^\num$. Пусть
$$
F:=F_{p^*,p',p'';q^*,q',q'';r^*,r',r''}(M,\d^+M,\d^-M)
$$
-- пространство функций Морса $f\in F^\num$ на поверхности
$(M,\d^+M,\d^-M)$,
таких что $\N^*_{f,\lambda}=\N^*_{f_*,\lambda}$ для любого
$\lambda=0,1,2$. Пространство $F$ мы наделим $C^\infty$-топологией,
см.~\cite[\S4]{kp1}, и назовем его {\em обобщенным пространством
функций Морса на по\-верх\-нос\-ти $(M,\d^+M,\d^-M)$}.
Обозначим через $F^1\subset F$ подпространство в $F$, состоящее из
таких функций Морса $f\in F$, что все локальные минимумы равны
$f(\d^-M)=-1$, а все локальные максимумы равны $f(\d^+M)=1$.
 \end{Def}

Из теоремы С.В.\ Матвеева (см.~\cite{Kmsb}) и ее обобщения в~\cite
{Kmsb} следует, что любое обобщенное пространство
$F=F_{p^*,p',p'';0,\hat q,q'';r^*,r',r''}(M,\d^+M,\d^-M)$ функций
Морса без закрепленных седловых точек (т.е.\ при $q^*=0$) линейно
связно.

\begin{Not} \label{not:R:G0*}
(A) Обозначим через $\N_\lam:=\N^*_{f,\lam}$ множество фиксированных
критических точек индекса $\lam$ (совпадающее для разных функций
$f\in F$), $\lam=0,1,2$, положим $\N:=\N_0\cup\N_1\cup\N_2$. Пусть
$\D=\Diff^+(M,\d^+M,\d^-M,\N_0,\N_1,\N_2)$ --- группа сохраняющих ориентацию
диффеоморфизмов поверхности $M$, переводящих каждое множество
$\d^+M,\d^-M$, $\N_\lam$ в себя, $\lam=0,1,2$.
Пусть $\D^0=\Diff^0(M,\N)$ --- множество
всех диффео\-мор\-физ\-мов $h\in\D$, гомотопных $\id_M$ в классе
гомеоморфизмов пары $(M,\N)$. Пространства
$\D^0\subset\D$ наделим $C^\infty$-топологией,
см.~\cite[\S4(б)]{kp1}.

(B) Обозначим через $\bar M$ замкнутую поверхность, полученную из
поверхности $M$ стягиванием в точку каждой граничной окружности.
Обозначим через $\T\subset\D$ группу (называемую {\em группой
Торелли}), состоящую из всех диффеоморфизмов $h\in\D$, переводящих в
себя каждую компоненту края $M$, и таких что индуцированный
гомеоморфизм $\bar h\:\bar M\to\bar M$ индуцирует тождественный
автоморфизм группы гомо\-ло\-гий $H_1(\bar M)$. Имеем
$\D^0\subset\T$.
 \end{Not}

Из результатов~\cite {EE,EE0} следует, что имеется гомотопическая эквивалентность
 \begin{equation} \label {eq:EE}
 \D^0 \aapprox R_{\D^0},
 \end{equation}
где $R_{\D^0}$ --- одно из четырех многообразий, определяемое парой
$(M,|\N|)$, а именно: $SO(3)=\RR P^3$ (при $M=S^2$, $\N=\emptyset$),
$SO(2)=S^1$ (при $0\le\chi(M)-|\N|\le1$ и $d^++d^-+|\N|>0$),
$T^2=S^1\times S^1$ (при $M=T^2$, $\N=\emptyset$) и точка (при
$\chi(M)<|\N|$) (см., например,~\cite{S,EE}). В частности, $\D^0$
линейно связно. Кроме того,
 \begin{equation} \label {eq:EE1}
 \D^0=\T \quad \iff \quad |\N|\le\chi(M)+1.
 \end{equation}
Импликация ``$\Leftarrow$'' в (\ref {eq:EE1}) следует из~\cite
{EE,EE0}, а импликация ``$\Rightarrow$'' следует из того, что в
случае $|\N|\ge\chi(M)+2$ существует диффеоморфизм $h\in\T$
(скручивание Дэна~\cite{Dehn} вокруг разбивающей окружности),
негомотопный $\id_M$ в пространстве гомео\-мор\-физ\-мов пары
$(M,\N)$, см.~\cite[лемма 2.1(1)]{BLMC} или \cite{Pdiplom}.

\begin{Def} \label{def:equiv}
{\rm(A)} Функции Морса $f,g\in F$ назовем {\em эквивалентными}
($f\sim g$), если найдутся такие диффеоморфизмы $h_1\in\D$ и
$h_2\in\Diff^+(\RR)$, что $f=h_2\circ g\circ h_1$ и $h_1$ сохраняет
нумерацию критических точек. Пусть $[f]$ -- класс эквивалентности.

{\rm(B)} Две функции Морса $f$ и $g$ назовем {\em изотопными}, если
они эквивалентны и $h_1\in\D^0$ (т.е.\ $h_1$ изотопен
тождественному), и обозначаем $f\sim_\isot g$. Множество всех функций
из $F^1$, изотопных функции $f$, обозначим через $[f]_\isot$.
\end{Def}

Классификация функций Морса из $F$ с точностью до (послойной)
экви\-ва\-лент\-нос\-ти
изучена в~\cite[гл.\,2, теоремы~4 и~8]{BFbook}, с точностью до
изотопности в~\cite[лемма~1 и теорема~2]{KP2} и~\cite[утверждение~1.1
и~\S3]{K}.

\begin{Not} \label{not:Gf:or}
Для любой функции Морса $f\in F$ рассмотрим граф $G_f$ в поверхности
$\Int(M)$, полученный из графа $f^{-1}(f(\N_{f,1}))$ выкидыванием
всех компо\-нент связности, не содержащих седловых критических точек
(см.\ определение \ref {def:Morse}). Этот граф имеет $q$ вершин
(являющихся седловыми точками $\b\in\N_{f,1}$), степени всех вершин
равны $4$, а значит в графе $2q$ ребер. Если поверхность $M$
ориентирована, то на ребрах графа $G_f$ имеется естественная
ориентация, такая, что в любой внут\-рен\-ней точке ребра репер,
составленный из положительно ориентированного ка\-са\-тель\-но\-го
вектора к ребру и вектора $\grad f$ (по отношению к какой-нибудь
фиксированной римановой метрике), задает положительную ориентацию
поверхности. Аналогично вводится ориентация на любой связной
компоненте линии уровня $f^{-1}(a)$ функции $f$, не содержащей
критическую точку, $a\in\RR$. Обозначим через $s(f):=|f(\N_{f,1})|$
количество седловых критических значений функции $f$.
\end{Not}

Пусть
$$
\MM=\MM_{p^*+d^-,p',p'';q^*,q',q'';r^*+d^+,r',r''}
$$
-- $3q$-мерное многообразие, содержащее комплекс
$\KK=\KK_{p^*+d^-,p',p'';q^*,q',q'';r^*+d^+,r',r''}$ оснащенных
функций Морса (см.\ \S\ref {sec:KK0} или~\cite[\S4]{KP3}). Пусть
$\MM_{[f]_\isot}\subset\MM_{\succeq[f]_\isot}$ --
$(s([f])+2q)$-мерный страт и его специальная окрестность в $\MM$,
отвечающие классу изотопности $[f]_\isot$ (см.\ утверждение \ref
{sta:MM1} или~\cite[\S4]{KP3}). Из \cite {KP3} нетрудно выводится,
что страт $\MM_{[f]_\isot}$ имеет своим сильным деформационным
ретрактом пространство орбит $(S^1)^{d([f])}/\Gamma_{[f]}$
соответствующего тора $(S^1)^{d([f])}$ по свободному действию
конечной группы $\Gamma_{[f]}$ допустимыми автоморфизмами тора, см.\
\cite[\S2]{KP3}.

\begin{Thm}  \label{thm:KP4add}
Пусть $M$ --- связная компактная ориентируемая поверхность с
разбиением края $\d M=\d^+M\sqcup\d^-M$ на положительные и
отрицательные окруж\-но\-сти. Рассмотрим обобщенные пространства
$$
F=F_{p^*,p',p'';q^*,q',q'';r^*,r',r''}(M,\d^+M,\d^-M), \qquad
F^1\subset F
$$
функций Морса на поверхности $(M,\d^+M,\d^-M)$, {\rm см.\
определение~$\ref {def:Morse}$} (т.е.\ у функ\-ций $f\in F$ некоторые
из критических точек могут быть отмечены, а некоторые закреплены).
Пусть $\FF^1\subset\FF$ --- соответствующие пространства оснащенных
функ\-ций Морса {\rm(см.\ определение~$\ref {def:framed}$)}.
Предположим, что
 \begin{equation} \label {eq:main}
\hat p+\hat q+\hat r >\chi(M)
 \end{equation}
(т.е.\ количество отмеченных критических точек превосходит
$\chi(M)$). Тогда:

{\rm (A)} Имеются гомотопические эквивалентности и гомеоморфизм
 $$
 F \sim F^1 \sim \FF \sim \FF^1
 \approx \D^0\times\MM \qquad  (\aapprox R_{\D^0}\times\MM),
 $$
где $R_{\D^0}$ -- одно из многообразий $\RR P^3$, $S^1$, $S^1\times S^1$
и точка, см.\ {\rm(\ref {eq:EE})}.

{\rm (B)} Для любой функции Морса $f\in F^1$ имеются гомотопические эквивалентности и гомеоморфизм
 $$
 [f]_\isot\sim \Forg_1^{-1}([f]_{\isot}) \approx \D^0\times \MM_{[f]_\isot}
 \aapprox \D^0\times ((S^1)^{d}/\Gamma_{[f]})
 \qquad (\aapprox R_{\D^0}\times ((S^1)^{d}/\Gamma_{[f]})),
 $$
где $\Forg_1\:\FF^1\to F^1$ -- забывающее отображение, $\MM_{[f]_\isot}\subset\MM$ и
$(S^1)^d=(S^1)^{d([f])}$ -- соответствующие $(s([f])+2q)$-мерное подмногообразие и тор.
 \end{Thm}

Пусть $\Bbbk$ -- поле (например, $\RR,\QQ$ или $\ZZ_p$). Для
топологического пространства $X$ рассмотрим его числа Бетти
$\beta_j(X):=\dim_\Bbbk H_j(X;\Bbbk)$ и полином Пуанкаре
$P(X,t):=\sum_{j=0}^\infty t^j\beta_j(X)$. Следующее утверждение
(аналогичное \cite [следствие 2.7]{KP3}) выводится из теоремы \ref
{thm:KP4add} и определения \ref {def:MM} стратифицированного
многообразия $\MM$ стан\-дарт\-ными методами теории Морса (см.,
например,~\cite[\S45]{FF}), при помощи фильтрации $\emptyset=\MM_{\ge
q+1}\subset\MM_{\ge q}\subset\ldots\subset\MM_{\ge1}=\MM$ и рассмотрения
индуцированной стратификации специальной окрестности
$\MM_{\succeq[f]_\isot}$ каждого страта $\MM_{[f]_\isot}$ в $\MM$
(см.\ утверждение \ref {sta:MM1}), где $\MM_{\ge s}:=\cup_{s([f])\ge
s}\MM_{\succeq[f]_\isot}=\cup_{s([f])\ge s}\MM_{[f]_\isot}$.

\begin{Cor} \label {cor:ineq}
{\rm(A)} Если количество $\hat p+\hat q+\hat r$ отмеченных
критических точек превосходит $\chi(M)$, то $\beta_j(F)=0$ при любом
$j\ge3q+2$.

{\rm(B)} Пусть $\bar M=S^2$ {\rm(см.\ обозначение \ref
{not:R:G0*}(B))}, $p^*+q^*+r^*\le\chi(M)+1\le\hat p+\hat q+\hat r$.
Тогда $\D=\D^0$; стратифицированное $3q$-мерное многообразие
$\IMM:=\MM$ состоит из конечного числа стратов
$\IMM_{[f]}:=\MM_{[f]_\isot}$;
имеется гомотопическая экви\-ва\-лент\-ность $F\sim R\times\IMM$, где
$R$ -- одно из многообразий $\RR P^3$, $S^1$, $S^1$ и точка в
зависимости от значения $\chi(M)-(p^*+q^*+r^*)=2,1,0,-1$; числа Бетти
$\beta_j=\beta_j(\IMM)$ многообразия $\IMM$ удовлетворяют
неравенствам Морса-Смейла:
 $$
 \beta_j-\beta_{j-1}+\beta_{j-2}-\beta_{j-3}+\ldots\le
 q_j-q_{j-1}+q_{j-2}-q_{j-3}+\ldots, \qquad j\ge0,
 $$
где $Q(t)=\sum_{j=0}^\infty t^jq_j:=\sum_{[f]\in
F^1/\sim}t^{q-s(f)}P(\IMM_{[f]},t)$.
В частности, справедливы неравенства Морса:
 $$
 \chi(\IMM)=(-1)^{q-1}\left|\left\{[f]\in F^1/\sim \ \mid\
s(f)=1\right\}\right|, \qquad \beta_j\le q_j, \quad j\ge0.
 $$
\end{Cor}

\subsection{Комбинаторное построение многообразия $\MM$ согласно \cite {KP3}}\label {sec:KK0}

\begin{Not} \label {not:numeration}
(A) Аналогично определению \ref {def:Morse} и обозначению \ref
{not:R:G0*} обо\-зна\-чим через $\N_f := \N_{f,0} \cup \N_{f,1} \cup \N_{f,2}$,
$\hat\N_f:=\hat\N_{f,0}\cup \hat\N_{f,1}\cup \hat\N_{f,2}$
множество всех критических точек (соответственно всех отмеченных
критических точек) функции $f\in F$. Имеем включения
$\N\subseteq\hat\N_{f}\subseteq\N_{f}$ и
$\N_\lam\subseteq\hat\N_{f,\lam}\subseteq\N_{f,\lam}$ множеств
фиксированных критических точек, отмеченных критических точек и всех
критических точек (соответственно индекса $\lambda$) функции $f$,
$\lam=0,1,2$.

(B) Пусть $\sigma,\tau\subset X$ -- два непересекающихся подмножества
топологического про\-стран\-ст\-ва $X$ (например, две открытые клетки
клеточного комплекса). Будем гово\-рить, что $\sigma$ {\it примыкает
к} $\tau$ и писать $\tau\prec\sigma$ (и $\bar\tau\prec\bar\sigma$),
если $\tau\subset\d\sigma:=\bar\sigma\setminus\sigma$. Пишем
$\tau\preceq\sigma$, если $\tau\prec\sigma$ или $\tau=\sigma$.
\end{Not}

{\it Шаг 1.}
Пусть $J=(J_1,\dots,J_s)$ -- упорядоченное разбиение множества
$\{1,\dots,q\}$ на $s$ непустых подмножеств $J_1,\dots,J_s$ (т.е.\
$\{1,\dots,q\}=J_1\sqcup \ldots \sqcup J_s$), где $1\le s\le q$.
Определим числа $0=r_0<r_1<\ldots<r_{s-1}<r_s=q$ и перестановку
$\pi\in\Sigma_q$ условиями
 \begin{equation}\label{eq:J}
J_1=\{\pi_1,\dots,\pi_{r_1}\}, \ J_2=\{\pi_{r_1+1},\dots,\pi_{r_2}\},
\ \ldots, \ J_s=\{\pi_{r_{s-1}+1},\dots,\pi_{r_s}\},
 \end{equation}
$\pi_1<\ldots<\pi_{r_1}$, $\pi_{r_1+1}<\ldots<\pi_{r_2}, \ldots,
\pi_{r_{s-1}+1}<\ldots<\pi_{r_s}$.
Если разбиение $\hat J$ получается из разбиения $J=(J_1,\dots,J_s)$
путем измельчения (т.е.\ разбиения неко\-то\-рых множеств $J_k$ на
несколько подмножеств), будем писать $\hat J\prec J$.

{\it Шаг 2.} Для каждой функции Морса $f\in F$ рассмотрим множество
$\N_{f,1}=:\{y_j\}_{j=1}^q\approx\{1,\dots,q\}$ ее седловых
критических точек (см.\ обозначение~\ref {not:numeration}) и
евкли\-дово векторное пространство $0$-коцепей
 \begin{equation} \label{eq:H0f}
 H_f^0:=C^0(\N_{f,1};\RR)=\RR^{\N_{f,1}}\cong\RR^q
 \end{equation}
со стандартной евклидовой метрикой. Рассмотрим в пространстве $H^0_f$
внутренность куба: $(-1;1)^{\N_{f,1}}\approx(-1;1)^q\subset\RR^q$.
Рассмотрим ``вычисляющую'' 0-коцепь
 $$
 \barc=\barc(f):=f|_{\N_{f,1}}=(c_1,\dots,c_q)\in(-1;1)^{\N_{f,1}}\subset H_f^0,
 $$
т.е.\ функцию $\barc\:\N_{f,1}\to\RR$, сопоставляющую любой седловой
точке $y_j\in \N_{f,1}$ значение $c_j:=f(y_j)$ функции $f$ в этой
точке, $1\le j\le q$. Сопоставим 0-коцепи $\barc=(c_1,\dots,c_q)$
число $s(\barc):=|\{c_1,\dots,c_q\}|$ различных седловых значений и
упорядоченное разбиение $J=J(\barc)=(J_1,\dots,J_s)$ множества
седловых точек $\N_{f,1}\approx\{1,\dots,q\}$, определяемое
свойствами~(\ref {eq:J}) и $c_{\pi_1}=\ldots=c_{\pi_{r_1}} <
 c_{\pi_{r_1+1}}=\ldots=c_{\pi_{r_2}} < \ldots <
c_{\pi_{r_{s-1}+1}}=\ldots=c_{\pi_{r_s}}$. (То есть, $J$ -- это
отношение частичного порядка на множестве $\N_{f,1}$ седловых
критических точек функции $f$ значениями функции $f|_{\N_{f,1}}$.)

В каждом
классе изо\-топ\-нос\-ти $[f]_\isot\in F^1/\sim_\isot$ отметим ровно
одну функцию Морса $f$ этого класса.
Сопоставим классу изотопности $[f]_\isot$
и любому разбиению $J$ соответ\-ству\-ю\-щие страт и звездообразную
область в кубе $(-1;1)^{\N_{f,1}}$:
 $$
S_{f,J} := \{\barc'\in(-1;1)^{\N_{f,1}}\mid J(\barc')=J\}, \quad
S_{f,\preceq J} := \{\barc'\in(-1;1)^{\N_{f,1}}\mid J(\barc')\preceq
J\},
$$
$$
 S_{[f]_\isot} =S_{f}:=S_{f,J(\barc(f))}, \quad
 S_{\succeq[f]_\isot}
 :=S_{f,\preceq J(\barc(f))}.
 $$
Рассмотрим также двойственные друг другу векторные пространства
относи\-тель\-ных 1-гомологий и относительных 1-когомологий над полем
$\RR$:
 \begin{equation} \label {eq:H1f}
 \begin{array}{l}
 H_{f,1}:=H_1(M\setminus(\N_{f,0}\cup \N_{f,2}),\N_{f,1};\RR)\cong\RR^{2q}, \\
 \phantom{0} H^1_f:=H^1(M\setminus(\N_{f,0}\cup \N_{f,2}),\N_{f,1};\RR) \cong \Hom_\RR(H_{f,1},\RR) \cong \RR^{2q}.
 \end{array}
 \end{equation}
Рассмотрим ориентированный граф $G_f\subset M\setminus(\N_{f,0}\cup
\N_{f,2})$, см.\ обозначение~\ref {not:Gf:or}. Он имеет
$2q$ ребер, которые обозначим
$e_1,\dots,e_{2q}$. Обозначим относительный гомологический класс
ориентированного ребра $e_i$ через $[e_i]\in H_{f,1}$, $1\le i\le2q$.
Определим в векторном пространстве $H^1_f\cong\RR^{2q}$ выпуклое
подмножество
 \begin{equation} \label {eq:U0}
 U_{[f]_\isot}^{\infty}=U_f^{\infty}:=\left\{u\in H^1_f \ \left|\ u([e_i]) >0,\ 1\le i\le 2q \right. \right\}.
 \end{equation}

Через $\stab_{\D^0}g$ обозначим группу изотропии элемента $g\in F^1$
относительно естест\-вен\-но\-го правого действия группы $\D^0$ на
$F^1$, а через $(\stab_{\D^0}f)^0$ обозначим ее под\-груп\-пу,
состоящую из всех диффеоморфизмов поверхности $M$, сохраняющих
функ\-цию $f$ и гомотопных $\id_M$ в классе гомеоморфизмов $M$,
сохраняющих функцию $f$. Рассмотрим покомпонентное правое действие
дискретной группы
\begin{equation}\label{eq:Gamma:f}
 \widetilde\Gamma_{[f]_\isot}=\widetilde\Gamma_f:=(\stab_{\D^0}f)/(\stab_{\D^0}f)^0
 \end{equation}
на прямом произведении $S_{\succeq[f]_\isot}\times
U_{[f]_\isot}^\infty$ индуцированными автоморфизмами про\-странств
(\ref {eq:H0f}), (\ref {eq:H1f}). Согласно \cite[3.3]{KP3}, это
действие свободно и дискретно, а про\-стран\-ства орбит
 \begin{equation} \label{eq:DD}
 \MM_{\succeq[f]_\isot}^\st
 := (S_{\succeq[f]_\isot}\times U_f^\infty)/\widetilde\Gamma_f, \qquad
 \MM_{[f]_\isot}^\st
 := (S_{[f]_\isot}\times U_f^\infty)/\widetilde\Gamma_f \subset\MM_{\succeq[f]_\isot}^\st
 \end{equation}
являются $3q$-мерным открытым многообразием и его $(s(f)+2q)$-мерным
под\-мно\-го\-о\-бра\-зием соответственно.

{\it Шаг 3.} Изучим взаимосвязь $3q$-мерных многообразий
$\MM_{\succeq[f]_\isot}^\st$, $\MM_{\succeq[g]_\isot}^\st$ для
примыкающих классов изотопности $[f]_\isot\prec[g]_\isot$ (см.\
обозначение~\ref {not:numeration}(B)). Пусть $f\in F^1$ -- отмеченная
функция своего класса изотопности, и пусть функция $\tilde f\in F^1$
получена малым возмущением функции $f\in F^1$, причем $\N_{\tilde
f}=\N_f$. Обозначим через $g$ отмеченную функцию класса изотопности
$[\tilde f]_\isot$, $J':=J(\barc(\tilde f))$, и через
$h_{[f]_\isot,J'}=h_{f,J'}:=h_{\tilde f,g}\in\D^0$ диффеоморфизм,
переводящий линии уровня функции $g$ в линии уровня функции $\tilde
f$ с сохранением направления роста (он существует ввиду изотопности
функций $\tilde f,g$). Согласно \cite[утверждение~1.1 и~\S3]{K} или
\cite{KP3}, выполнено $J'\preceq J:=J(\barc(f))$ и имеется сюръекция
$\delta[f]_\isot$ множества всех упорядоченных разбиений $J'\preceq
J$ на множество всех классов изотопности $[g]_\isot\succeq[f]_\isot$
(см.\ обозна\-че\-ние~\ref {not:numeration}(B)), такая что
 \begin{equation} \label {eq:incid}
 \delta{[f]_\isot}\:J'=J(\barc(\tilde f))=J((h_{f,J'}^{-1})^{*0}(\barc(g)))\mapsto\delta_{J'}[f]_\isot:=[\tilde f]_\isot.
 \end{equation}
Хотя сопоставление $([f]_\isot,J')\mapsto h_{f,J'}$ не является
однозначным (т.е.\ диффео\-мор\-физм $h_{f,J'}$ зависит, вообще
говоря, от возмущенной функции $\tilde f$, такой что $\N_{\tilde
f}=\N_f$ и $J(\barc(\tilde f))=J'$), но в силу~\cite[лемма 1]{KP2}
смежный класс $h_{f,J'}(\stab_{\D^0}g)(\Diff^0(M,\N_g))$ определен
корректно, где через $\Diff^0(M,\N_{g})\subset\D^0$ обозначена группа
диффеоморфизмов пары $(M,\N_{g})$, гомотопных $\id_M$ в классе
гомеоморфизмов пары.

Рассмотрим индуцированные изоморфизмы векторных пространств:
\begin{equation}\label{eq:embedU}
h_{f,J'}^{*0}\:H^0_f\to H^0_g, \qquad h_{f,J'}^*\:H^1_f\to H^1_g,
 \end{equation}
см.\ (\ref {eq:H0f}), (\ref {eq:H1f}). Рассмотрим в
$(-1;1)^{\N_{f,1}}$ открытые $\widetilde\Gamma_f$-инвариантные
подмножества
 $$
 \d_{g} S_{\succeq[f]_\isot}=
 \d_{[g]_\isot} S_{\succeq[f]_\isot}
 :=
 \bigcup_{J'\in(\delta[f]_\isot)^{-1}([g]_\isot) } S_{f,\preceq J'} \subset S_{\succeq[f]_\isot}.
 $$
Согласно \cite [\S3]{KP3}, прямое произведение изоморфизмов в
(\ref {eq:embedU}) индуцирует корректно определенное вложение
$3q$-мерных открытых многообразий:
 $$
  \chi_{[f]_\isot,[g]_\isot}= \chi_{f,g}
 \: \d_{[g]_\isot}\MM_{\succeq[f]_\isot}^\st:=
 \left(\left(\d_gS_{\succeq[f]_\isot}\right)\times U_f^\infty\right)/\tilde\Gamma_f
 \hookrightarrow
\MM_{\succeq[g]_\isot}^\st,
 $$
 $$
 \tilde\Gamma_f(\barc,u) \mapsto \tilde\Gamma_g(h_{f,J'}^{*0}(\barc),h_{f,J'}^*(u)),
 \qquad (\barc,u)\in S_{f,\preceq J'}\times U_f^\infty,
 $$
где $\delta_{J'}[f]_\isot=[g]_\isot$ (см.~(\ref {eq:DD}), (\ref
{eq:incid})).

{\it Шаг 4.}
Предположим, что отмеченные функции $f$ всех классов изотопности
$[f]_\isot$ имеют одно и то же
множество критических точек $\N_{f,\lam}=\N_{f_*,\lam}$ с учетом меток, $\lam=0,1,2$
(см.\ определение~\ref {def:Morse}(B,C) и шаг 2).
Рассмотрим топологическое пространство
$$
 (F^1/\sim_\isot)^\diskr \ \times (-1;1)^{\N_{f_*,1}}\times H^1_{f_*}
 \ \approx \
 (F^1/\sim_\isot)^\diskr \ \times (-1;1)^q\times\RR^{2q},
$$
где $(F^1/\sim_\isot)^\diskr:=F^1/\sim_\isot$ с дискретной
топологией, и его подпространства
$$
 \XX_\succeq:=\bigcup\limits_{[f]_\isot\in F^1/\sim_\isot} \{[f]_\isot\}\times S_{\succeq[f]_\isot}
 \times U_{[f]_\isot}^{\infty},
$$
$$
 \XX:=\bigcup\limits_{[f]_\isot\in F^1/\sim_\isot} \{[f]_\isot\}\times S_{[f]_\isot}
 \times U_{[f]_\isot}^{\infty} \subset \XX_\succeq.
$$

\begin{Def} [многообразие $\MM$] \label {def:MM}
Пусть $\MM:=(\XX_\succeq/\sim)/\sim_\glue$ -- про\-стран\-ство с
фактортопологией, где отношения эквивалентности $\sim$, $\sim_\glue$
на множествах $\XX_\succeq$, $\YY_\succeq:=\XX_\succeq/\sim$
порождены следующими отношениями:

(отношение $\sim$ на $\XX_\succeq$) для каждого класса изотопности $[f]_\isot$ рассмотрим проек\-цию
$\{[f]_\isot\}\times S_{\succeq[f]_\isot}\times U_{[f]_\isot}^\infty
 \to
\{[f]_\isot\}\times\MM^\st_{\succeq[f]_\isot}=:\uups_{\succeq[f]_\isot}$
(см.\ (\ref {eq:DD})),
и назовем точки множества $\{[f]_\isot\}\times S_{\succeq[f]_\isot}
\times U_{[f]_\isot}^\infty$ {\it $\sim$-эквивалентными}, если их
образы при этой проекции совпадают; обозначим
$\uups_{[f]_\isot}:=\{[f]_\isot\}\times\MM^\st_{[f]_\isot}\subset\uups_{\succeq[f]_\isot}\subset\YY_\succeq$;

(отношение $\sim_\glue$ на $\YY_\succeq$; отображения инцидентности)
для каждой пары при\-мы\-ка\-ю\-щих классов $[f]_\isot\prec[g]_\isot$
(см.\ обозначение \ref {not:numeration}) рассмотрим вложение
соот\-вет\-ству\-ю\-щих $3q$-мерных открытых многообразий, называемое
{\it отображением инци\-дент\-нос\-ти} этой пары:
$\chi_{[f]_\isot,[g]_\isot}\:\d_{[g]_\isot}\MM^\st_{\succeq[f]_\isot}\hookrightarrow\MM^\st_{\succeq[g]_\isot}$;
рассмотрим инду\-ци\-ро\-ван\-ное вложение
$\d_{[g]_\isot}\uups_{\succeq[f]_\isot}:=\{[f]_\isot\}\times(\d_{[g]_\isot}\MM^\st_{\succeq[f]_\isot})\hookrightarrow\uups_{\succeq[g]_\isot}$
(которое тоже обозначим $\chi_{[f]_\isot,[g]_\isot}$); назовем любую
точку множества
$\d_{[g]_\isot}\uups_{\succeq[f]_\isot}\subset\YY_\succeq$ и ее образ в
$\uups_{\succeq[g]_\isot}\subset\YY_\succeq$ при данном вложении {\it
$\sim_\glue$-эквивалентными}.
\end{Def}

Пусть $p_X\:\XX_\succeq\to\MM$, $p_Y\:\YY_\succeq\to\MM$ --
канонические проекции. Положим $\YY:=\XX/\sim\subset\YY_\succeq$,
$\MM_{\succeq[f]_\isot}:=p_Y(\uups_{\succeq[f]_\isot})$, назовем
$\MM_{[f]_\isot}:=p_Y(\uups_{[f]_\isot})$ {\it стратом} в $\MM$.

Так как $\XX_\succeq$ -- гладкое открытое $3q$-мерное многообразие с
естественной плоской аффинной связностью, и группы
$\widetilde\Gamma_{[f]_\isot}$ действуют на нем с сохранением
связности, то $\YY_\succeq$ тоже является гладким открытым
$3q$-мерным многообразием с плоской аф\-фин\-ной связностью, причем
пересечение любой его связной компоненты $\uups_{\succeq[f]_\isot}$ с
подмножеством $\YY\subset\YY_\succeq$ является плоским
$(s([f])+2q)$-мерным подмногообразием.

\begin{Sta} [{\rm \cite [теорема 4.3]{KP3}}] \label {sta:MM1}
Пространство $\MM:=\YY_\succeq/\sim_\glue$ обладает структурой
гладкого $3q$-мерного многообразия и плоской аффинной связностью,
гладкой относительно этой структуры. Для каждого класса изотопности
$[f]_\isot$ отображение
$p_Y|_{\uups_{\succeq[f]_\isot}}\:\uups_{\succeq[f]_\isot}\to\MM$
является гладким регулярным вложением гладких $3q$-мерных
многообразий, сохраняющим аффинную связность, а потому подмножество
$\MM_{\succeq[f]_\isot}=p_Y(\uups_{\succeq[f]_\isot})\subset\MM$
открыто.
Отображение $p_Y|_\YY\:\YY\to\MM$ биективно. В частности, страты
$\MM_{[f]_\isot}\subset\MM$ попарно не пересекаются, являются
плоскими $(s([f])+2q)$-мерными подмногообразиями и покрывают все
$\MM$. Дискретная группа $\D/\D^0$ и группа $\Diff^+[-1;1]$ действуют
на $\MM$ справа и слева (соответственно) диффеоморфизмами,
сохраняющими аффинную связность, стратификацию и систему открытых
подмножеств $\MM_{\succeq[f]_\isot}\subset\MM$ (называемых
специальными окрестностями стратов).
\end{Sta}

\subsection{Гомеоморфизм между универсальным пространством модулей $\FF^1/\D^0$
оснащенных функций Морса и многообразием $\MM$} \label{sec:MM}

\begin{Def} [{\rm\cite[\S9]{kp1}}] \label{def:framed}
{\em Оснащенной функцией Морса на ориентированной поверхности \\
$(M,\d^+M,\d^-M)$} назовем пару $(f,\alpha)$, где $f\in F$ ---
функция Морса на $(M,\d^+M,\d^-M)$, $\alpha$ --- замкнутая 1-форма на
$M\setminus(\N_{f,0}\cup\N_{f,2})$, такие что $2$-форма
$df\wedge\alpha$ не имеет нулей в $M\setminus\N_f$ и задает
положительную ориентацию, и в окрестности любой
критической точки $x\in\N_{f}$ существуют
локальные координаты $u,v$, в которых либо $f=u^2-v^2+f(x)$,
$\alpha=d(2uv)$, либо $f=\varkappa_{f,x}(u^2+v^2)+f(x)$,
$\alpha=\varkappa_{f,x}\frac{udv-vdu}{u^2+v^2}$, где
$\varkappa_{f,x}=\const\ne0$. Обозначим через $\FF^1$ пространство
оснащенных функций Морса $(f,\alpha)$, таких что $f\in F^1$. Снабдим
его $C^\infty$-топологией (см.~\cite[\S4]{kp1}).
\end{Def}

Сформулируем без доказательства техническую лемму.

\begin{Lem} \label {lem:1}
Для любой функции Морса $f\in F^1$ существует гладкое
$3q$-пара\-мет\-ри\-чес\-кое семейство оснащенных функций Морса
$(f_{\barc'},\alpha_{f,u})\in\FF^1$ с параметрами $(\barc',u)\in
S_{f,\preceq J(\barc(f))}\times U_f^{\infty}$, такое что
$f_{\barc(f)}=f$, $\N_{f_{\barc'}}=\N_f$, $\barc(f_{\barc'})=\barc'$, $[\alpha_{f,u}]=u\in H^1_f$.
\end{Lem}

Определим ``вычисляющее'' отображение $\Ev\:\FF^1\to\MM$ формулой
\begin{equation} \label{eq:Ev}
\Ev(f,\alpha):=p_X([f]_\isot,h_{f,f_0}^{*0}(\barc(f)),h_{f,f_0}^*[\alpha]),
\qquad (f,\alpha)\in\FF^1,
\end{equation}
где $f_0\in F^1$ -- отмеченная функция Морса класса изотопности
$[f]_\isot$, $h_{f,f_0}\in\D^0$ -- какой-нибудь диффеоморфизм,
переводящий линии уровня функции $f_0$ в линии уровня функции $f$ с
сохранением направления роста и нумерации отмеченных
кри\-ти\-чес\-ких точек (см.~(\ref {eq:incid})), а
$h_{f,f_0}^{*0}\:H^0_f\to H^0_{f_0}$ и $h_{f,f_0}^*\:H^1_f\to
H^1_{f_0}$ -- индуцированные изоморфизмы (см.\ (\ref {eq:H0f}), (\ref
{eq:H1f}), (\ref {eq:incid})).

\begin{Sta} \label {sta:MM2}
Отображение $\Ev\:\FF^1\to\MM$ однозначно, $\D^0$-ин\-ва\-ри\-антно,
непрерывно и индуцирует $\D/\D^0$-эквивариантный гомеоморфизм
$\overline{\Ev}\:\FF^1/\D^0\to\MM$.
\end{Sta}

\begin{proof}
{\it Шаг 1.} Проверим однозначность отображения $\Ev$.
Образ \\ $p_X([f]_\isot,h_{f,f_0}^{*0}(\barc(f)),h_{f,f_0}^*[\alpha])\in\MM$
точки
$([f]_\isot,h_{f,f_0}^{*0}(\barc(f)),h_{f,f_0}^*[\alpha])\in\{[f]_\isot\}\times
S_{[f]_\isot}\times U_f^\infty\subset\XX$
не зависит от выбора диффеоморфизма $h_{f,f_0}$, так как для любого другого
такого диффеоморфизма $\tilde h_{f,f_0}$ в силу~\cite[лемма 1]{KP2} выполнено $\tilde
h_{f,f_0}^{-1}h_{f,f_0}\in(\stab_{\D^0}f_0)(\Diff^0(M,\N_{f_0}))$,
а действие группы
$(\stab_{\D^0}f_0)(\Diff^0(M,\N_{f_0}))$ на
$\MM_{\succeq[f_0]_\isot}^\st\approx\uups_{\succeq[f_0]_\isot}\subset\YY_\succeq$
тривиально (см.~(\ref {eq:Gamma:f}), (\ref {eq:DD})). Однозначность
$\Ev$ доказана.

{\it Шаг 2.} Докажем непрерывность отображения $\Ev$ в любой точке
$(f,\alpha)\in\FF^1$.
Для оснащенной функции Морса $(\tilde f,\tilde\alpha)\in\FF^1$, достаточно
близкой к $(f,\alpha)$,
рассмотрим близкий к $\id_M$ диффеоморфизм $h\in\D^0$, такой что
$h(\N_f)=\N_{\tilde f}$, и упорядоченное разбиение $J'\preceq
J(\barc(fh_{f,f_0}))$, такое что $J(\barc(\tilde fhh_{f,f_0}))=J'$, откуда
$\delta_{J'}[f]_\isot=[\tilde f]_\isot$, см.~(\ref {eq:incid}).
Пусть $g$ -- отмеченная функция класса изотопности $[\tilde f]_\isot$. Тогда
\begin{equation} \label {eq:*.}
 \Ev(\tilde f,\tilde\alpha)
 = p_X([g]_\isot,h^{*0}_{\tilde fhh_{f,f_0},g}(\barc(\tilde fhh_{f,f_0})),h^*_{\tilde fhh_{f,f_0},g}[h_{f,f_0}^*h^*\tilde\alpha])
\end{equation}
$$
 =p_X([g]_\isot,h^{*0}_{f_0,J'}(\barc(\tilde fhh_{f,f_0})),h^*_{f_0,J'}[h_{f,f_0}^*h^*\tilde\alpha])
 =p_X([f]_\isot,h_{f,f_0}^{*0}(\barc(\tilde fh)),h_{f,f_0}^*[h^*\tilde\alpha]),
$$
где последнее равенство следует из того, что отображение
инцидентности $\chi_{f_0,g}$ ин\-ду\-ци\-ро\-вано отображением
$h_{f_0,J'}\in h_{\tilde
fhh_{f,f_0},g}(\stab_{\D^0}g)(\Diff^0(M,\N_g))$. Из $C^0$-близости
$h$ к $\id_M$ следует, что 0-коцепь $\barc(\tilde fh)\in H^0_f$
близка к $\barc(f)$ (ввиду $C^2$-близости функции $\tilde f$ к $f$),
а класс относительных 1-когомологий $[h^*\tilde\alpha]\in H^1_f$
близок к $[\alpha]$ (ввиду $C^0$-бли\-зос\-ти 1-формы $\tilde\alpha$
к $\alpha$ вне малых окрестностей точек локальных минимумов и
мак\-си\-му\-мов функции $f$, см.\ определение топологии в
пространстве $\FF^1$ в~\cite[\S4.2]{kp1}). Поэтому точка $\Ev(\tilde
f,\tilde\alpha)$ близка к
$p_X([f]_\isot,h_{f,f_0}^{*0}(\barc(f)),h_{f,f_0}^*[\alpha])=\Ev(f,\alpha)$
ввиду непрерывности проекции
$p_X\:\XX_\succeq\to\MM:=(\XX_\succeq/\sim)/\sim_\glue$.
Непрерывность $\Ev$ доказана.

{\it Шаг 3.} По построению $\Ev$ является $\D^0$-инвариантным.
Индуцированное отобра\-же\-ние $\overline{\Ev}\:\FF^1/\D^0\to\MM$
непрерывно ввиду непрерывности отображения $\Ev$. Оно
$\D/\D^0$-эквивариантно по построению. Покажем, что $\overline{\Ev}$
биективно.

{\it Инъективность.} Пусть $\Ev(f,\alpha)=\Ev(f_1,\alpha_1)$. Ввиду
инъективности $p_Y|_{\YY}$ (см.\ утверждение \ref {sta:MM1})
выполнено $[f]_\isot=[f_1]_\isot$ и имеется диффеоморфизм $h_1\in\D^0$,
переводящий линии уровня функции $f$ в линии уровня функции $f_1$ с
сохранением направления роста и такой, что
$h_1^{*0}(\barc(f_1))=\barc(f)\in H^0_f$ и
$h_1^*[\alpha_1]\in\widetilde\Gamma_{[f]_\isot}[\alpha]$. Отсюда
$h_2^*h_1^*[\alpha_1]=[\alpha]$ для некоторого $h_2\in\stab_{\D^0}f$
(см.~(\ref {eq:Gamma:f})). Поэтому для
$(f_2,\alpha_2):=(h_1h_2)^*(f_1,\alpha_1)\in\FF^1$ выполнено
$G_{f_2}=G_f$, $\barc(f_2)=\barc(f)$, $[\alpha_2]=[\alpha]$.

Покажем, что существует (единственный) диффеоморфизм $h\in\D$,
переводящий в себя каждое ориентированное ребро графа $G_f$ и такой,
что $h^*(f_2,\alpha_2)=(f,\alpha)$.
В малых окрестностях $U_j,\tilde U_j$ каждой седловой точки $\b_j\in
\N_{f,1}$ в $M$ рассмотрим локальные координаты $u,v$ для
$(f,\alpha)|_{U_j}$ и $u_2,v_2$ для $(f_2,\alpha_2)|_{\tilde U_j}$
как в определении~\ref {def:framed}.
Без ограничения общности будем считать, что начальные отрезки вида
$\{0\le u=v\le\eps\}$ и $\{0\le u_2=v_2\le\eps_2\}$ ребер графа
$G_f$, выходящих из вершины $\b_j$, совпадают (в противном случае
заменим $(u_2,v_2)$ на $(-u_2,-v_2)$). Определим диффеоморфизм
$h|_{U_j'}\:U_j'\to\tilde U_j'$ в, быть может, меньшей окрестности
$U_j'\subset U_j$ условием $(u_2,v_2)\circ h|_{U_j'}=(u,v)|_{U_j'}$, где
$\tilde U_j':=h(U_j')\subset\tilde U_j$. Тогда
$h|_{U_j'}^*(f_2,\alpha_2)=(f,\alpha)|_{U_j'}$. Продолжим этот
диффеоморфизм на каждое ребро $e_\ell$ графа $G_f$ условием
$(h|_{e_\ell})^*(\alpha_2|_{e_\ell})=\alpha|_{e_\ell}$. Это возможно,
так как интегралы 1-форм $\alpha$ и $\alpha_2$ по ориентированному
ребру $e_\ell$ равны. Продолжим этот диффеоморфим в малую окрестность
$V_\ell$ куска $e_\ell\setminus(\cup_{j=1}^qU_j')$ этого ребра
условием $h|_{V_\ell}^*(f_2,\alpha_2)=(f,\alpha)|_{V_\ell}$, положим
$\tilde V_\ell:=h(V_\ell)$. На множестве $\N_{f,0}\cup \N_{f,2}$
точек локальных минимумов и максимумов определим $h|_{\N_{f,0}\cup
\N_{f,2}}:=\id_{\N_{f,0}\cup \N_{f,2}}$. Осталось продолжить
построенное отображение на открытое подмножество $M\setminus(G_f\cup
\N_{f,0}\cup \N_{f,2})\subset M$, являющееся дизъюнктным объединением
кусков, каждый из которых гомеоморфен открытому или полуоткрытому
цилиндру $S^1\times (0;1)$ и $S^1\times [0;1)$. Для каждого такого
куска $Z$ отображение $h$ уже построено на $\overline Z\setminus
Z\subset G_f\cup \N_{f,0}\cup \N_{f,2}$. Пусть точка $x\in\overline
Z\cap G_f$ не является критической, и пусть окружность
$\gamma_Z:=Z\cap (f^{-1}(\frac12(\inf f|_Z+\sup f|_Z)))$
ориентирована как в обозначении~\ref {not:Gf:or}. Для любой точки
$y\in Z$ рассмотрим гладкий путь $\gamma_{x,y}\:[0;1]\to\overline Z$
из $x$ в $y$, такой что $\gamma_{x,y}((0;1))\subset Z$. Положим
$$
A_{Z,\alpha}:=\oint_{\gamma_Z}\alpha>0, \qquad
g_{Z,x,\alpha}(y):=\int_{\gamma_{x,y}}\alpha\in\RR, \quad y\in Z.
$$
Так как $Z$ гомеоморфен открытому или полуоткрытому цилиндру и
1-форма $\alpha$ замкнута, то функция $g_{Z,x,\alpha}\mod
A_{Z,\alpha}\:Z\to\RR/A_{Z,\alpha}\ZZ$ корректно определена, т.е.\ не
зависит от выбора пути $\gamma_{x,y}$. По условию
$A_{Z,\alpha}=A_{Z,\alpha_2}$. Определим отображение $h|_Z$ условием
$(f_2,g_{Z,h(x),\alpha_2}\mod A_{Z,\alpha})\circ
h|_Z=(f,g_{Z,x,\alpha}\mod A_{Z,\alpha})$. Пусть $e_\ell\subset G_f$
-- ребро, содержащее точку $x$. Нетрудно доказывается непрерывность
отображения $h|_{Z\cup e_\ell}$. Отсюда, с учетом равенства
$[\alpha]=[\alpha_2]\in H^1_f$, следуют непрерывность и биективность
отображения $h|_{\overline Z}$. То, что $h|_{M\setminus(\N_{f,0}\cup
\N_{f,2})}$ является диффеоморфизмом, следует из того, что следующие
пары функций являются регулярными координатами: пара координат
$(u,v)$ в $U_j'$ (соответственно $(u_2,v_2)$ в $\tilde U_j'$); пара
функций $(f,g_\ell)$ в $V_\ell$ (соот\-вет\-ствен\-но $(f_2,\tilde
g_\ell)$ в $\tilde V_\ell$), где функции $g_\ell$ и $\tilde g_\ell$
определены условиями $dg_\ell=\alpha|_{V_\ell}$ и $d\tilde
g_\ell=\alpha_2|_{\tilde V_\ell}$; пара функций
$(f_2,g_{Z,h(x),\alpha_2}\mod A_{Z,\alpha})$ и $(f,g_{Z,x,\alpha}\mod
A_{Z,\alpha})$ в $Z$. По построению $(f,\alpha)=h^*(f_2,\alpha_2)$.
То, что $h$ является диффеоморфизмом в малой окрестности $W_\a$ любой
точки $\a\in \N_{f,0}\cup \N_{f,2}$ минимума или максимума,
доказывается с помощью полярных координат, отвечающих регулярным
координатам $u,v$ для $(f,\alpha)|_{W_\a}$ (см.\ определение~\ref
{def:framed}),
и аналогичных полярных координат для $(f_2,\alpha_2)|_{W_x}$.

Так как диффеоморфизм $h$ переводит в себя каждое ориентированное
ребро гра\-фа $G_f=G_{f_2}$, причем
$h^*[\alpha_2]=[\alpha]=[\alpha_2]\in U_f^\infty$, то
$h\in\Diff^0(M,\N_f)\subset\D^0$ согласно \cite [лемма 3.4]{KP3}.
По доказанному
$(f,\alpha)=h^*(f_2,\alpha_2)=(h_1h_2h)^*(f_1,\alpha_1)\in
\D^0(f_1,\alpha_1)\in\FF^1/\D^0$, и инъективность доказана.

{\it Сюръективность.} Отображение $p_X|_{\XX}\:\XX\to\MM$ сюръективно
ввиду сюръ\-ек\-тив\-нос\-ти отображений $\XX\to\YY=\XX/\sim$ и
$p_Y|_{\YY}\:\YY\to\MM$ (см.\ утверждение~\ref {sta:MM1}). Поэтому
достаточно показать, что для любой точки $([f]_\isot,\barc',u)\in\XX$
существует оснащенная функция Морса $(\tilde f,\tilde\alpha)\in\FF^1$, такая что
$p_X([f]_\isot,\barc',u)=\Ev(\tilde f,\tilde\alpha)$. Пусть
$f$ -- отмеченная функция своего класса изотопности. Из
включений $\barc'\in S_{f,J(\barc(f))}\subset S_{f,\preceq
J(\barc(f))}$, $u\in U^\infty_f$ и леммы~\ref{lem:1} получаем путь $(\tilde
f_t,\tilde\alpha):=(f_{t\barc'+(1-t)\barc(f)},\alpha_{f,u})\in\FF^1$
в пространстве $\FF^1$ оснащенных функций Морса, такой что
$\N_{\tilde f_t}=\N_f$, $\tilde f_0=f$, $\barc(\tilde f_1)=\barc'$ и
$[\tilde\alpha]=u$. Отсюда и из (\ref {eq:*.}) получаем требуемое
равенство
 $$
p_X([f]_\isot,\barc',u)=p_X([f]_\isot,\barc(f_{\barc'}),[\alpha_{f,u}])=\Ev(f_{\barc'},\alpha_{f,u}).
 $$

{\it Шаг 4.} Покажем, что непрерывная биекция
$\overline{\Ev}\:\FF^1/\D^0\to\MM$ является гоме\-о\-мор\-физ\-мом.
Осталось доказать, что ${\overline{\Ev}}^{-1}\:\MM\to\FF^1/\D^0$
непрерывно. Согласно лемме~\ref {lem:1} имеем непрерывное (ввиду
гладкости семейства) отображение
 \begin{equation} \label {eq:i4f}
\i_{f}\:
 S_{\succeq[f]_\isot}\times U_f^\infty\to\FF^1, \quad (\barc',u)\mapsto(f_{\barc'},\alpha_{f,u}),
 \end{equation}
для которого ввиду (\ref {eq:*.}) выполнено
\begin{equation}\label{eq:*1}
 \overline{\Ev}\circ q\circ\i_f(\barc',u) = \Ev\circ \i_{f}(\barc',u)=p_X([f]_\isot,\barc',u),
 \qquad (\barc',u)\in S_{\succeq[f]_\isot}\times U_f^\infty,
\end{equation}
где $q\:\FF^1\to\FF^1/\D^0$ -- проекция. Поэтому
$$
 {\overline{\Ev}}^{-1}\circ p_X|_{\{[f]_\isot\}\times S_{\succeq[f]_\isot}\times U_f^{\infty}}([f]_\isot,\barc',u)
 =q\circ \i_{f}(\barc',u), \quad (\barc',u)\in S_{\succeq[f]_\isot}\times U_f^\infty.
$$
Отсюда и из непрерывности отображения $q\circ\i_f$ следует
непрерывность отображения ${\overline{\Ev}}^{-1}$.
Действи\-тель\-но, по утверждению~\ref {sta:MM1} подмножество
$\MM_{\succeq[f]_\isot}$
открыто в $\MM$, а отображение \\ $p_X|_{\{[f]_\isot\}\times
S_{\succeq[f]_\isot}\times U_f^{\infty}}\:\{[f]_\isot\}\times
S_{\succeq[f]_\isot}\times U_f^{\infty}\to\MM_{\succeq[f]_\isot}$
есть композиция накрытия $\{[f]_\isot\}\times
S_{\succeq[f]_\isot}\times U_f^{\infty}\to\{[f]_\isot\}\times\MM_{\succeq[f]_\isot}^\st=\uups_{\succeq[f]_\isot}$
(см.\ (\ref {eq:DD})) и гомеоморфизма
$p_Y|_{\uups_{\succeq[f]_\isot}}\:\uups_{\succeq[f]_\isot}\to\MM_{\succeq[f]_\isot}$,
а потому оно локально является гомеоморфизмом. Утверждение~\ref {sta:MM2} доказано.
\end{proof}

\subsection{$\D^0$-эквивариантный гомеоморфизм $\FF^1\approx\D^0\times\MM$} \label {sec:FF:MM}

\begin{Not} \label {not:fix}
Предположим, что количество $|\hat\N_f\setminus\N|=\hat p+\hat q+\hat
r-(p^*+q^*+r^*)$ отмеченных, но не фиксированных, критических точек
любой функции $f\in F$ положительно. Фиксируем непустое подмножество
$\tilde\N_{f_*}\subseteq\hat\N_{f_*}\setminus\N$ и для
любой функции $f\in F$ обозначим через $\tilde\N_f\subseteq
\hat\N_f\setminus\N$ множество
ее
критических точек с теми же метками, что и точки множества
$\tilde\N_{f_*}$. Рассмотрим подпространство $F^*:=\{f\in
F\mid\tilde\N_f=\tilde\N_{f_*}\}\subset F$. Оно является обобщенным
пространством функций Морса на поверхности $(M,\d^+M,\d^-M)$ (см.\
определение~\ref {def:Morse}), каждая функция которого имеет ровно
$|\N|+|\tilde\N_f|=|\N^*|\in(|\N|,|\hat\N_{f_*}|]$ фиксированных критических точек, где
$\N^*:=\N\sqcup \tilde\N_{f_*}$ и $\N_\lam^*:=\N_\lam \cap \N^*$
суть множества всех фиксированных критических точек и фиксированных
критических точек индекса $\lam\in\{0,1,2\}$ соответственно,
$\N\subset\N^*\subseteq\hat\N_{f_*}$, см.\ обозначение \ref
{not:R:G0*}. Аналогично обозначению \ref {not:R:G0*} обозначим
 $$
\D^*:=\Diff^+(M,\d^+M,\d^-M,\N_0^*,\N_1^*,\N_2^*), \quad
(\D^*)^0:=\Diff^0(M,\N^*),
 $$
 $$
(\FF^*)^1:=\{(f,\alpha)\in\FF^1\mid f\in F^*\}.
 $$
В случае (\ref {eq:main}) рассмотрим для пространства $F^*$
соответствующее $3q$-мерное много\-обра\-зие $\MM^*$ (см.\
утверждение~\ref {sta:MM1}).
Пусть $\Ev^*\:(\FF^1)^*\to\MM^*$ -- вычисляющее отобра\-жение,
аналогичное вычисляющему отображению $\Ev\:\FF^1\to\MM$ (см.\
(\ref {eq:Ev})). По утверждению~\ref {sta:MM2} оно
индуцирует гомеоморфизм
$\overline{\Ev^*}\:(\FF^*)^1/(\D^*)^0\to\MM^*$.
\end{Not}

 \begin{Lem}\label{lem:KK*}
Для пространств $\FF^*\subset\FF$ обобщенных функций Морса {\rm(см.\
обозна\-че\-ние~\ref {not:fix})} отображения включения
$j\:\D^*\hookrightarrow\D$, $i\:(\FF^*)^1\hookrightarrow\FF^1$
индуцируют изо\-мор\-физм $\overline
j\:\D^*/(\D^*\cap\D^0)\stackrel{\cong}{\longrightarrow}\D/\D^0$ и
гомеоморфизм $\overline
i\:(\FF^*)^1/(\D^*\cap\D^0)\stackrel{\approx}{\longrightarrow}\FF^1/\D^0$.

Если $|\N^*|\le\chi(M)+1$, то $\D^*\cap\D^0=(\D^*)^0$, откуда имеются
изоморфизм
$\overline{j}\:\D^*/(\D^*)^0=\D^*/(\D^*\cap\D^0)\stackrel{\cong}{\longrightarrow}\D/\D^0$
групп классов отображений и гомео\-мор\-физм
$\overline{i}\:(\FF^*)^1/(\D^*)^0=(\FF^*)^1/(\D^*\cap\D^0)\stackrel{\approx}{\longrightarrow}\FF^1/\D^0$
универсальных пространств модулей, а в случае {\rm(\ref {eq:main})} также диффеоморфизм
$k:=\overline{\Ev}\circ\overline{i}\circ{\overline{\Ev^*}}^{-1}\:\MM^*\stackrel{\approx}{\longrightarrow}\MM$
$3q$-мерных многообразий (сохраняющий аффинную связность и
стратификацию).
\end{Lem}

\begin{proof}
Непосредственно проверяется, что $\overline j$ -- изоморфизм, а
$\overline i$ -- непрерывная биекция. Докажем непрерывность
$(\overline i)^{-1}$.
Для любой оснащенной функции Морса $(f_0,\alpha_0)\in\FF^1$
рассмотрим диффеоморфизм $h_0\in\D^0$, такой что
$h_0^*(f_0,\alpha_0)\in(\FF^*)^1$. Ввиду непрерывности отображения
$\FF^1\to F^1\to M^{|\tilde\N_{f_*}|}$, $(f,\alpha)\mapsto
f\mapsto\tilde\N_f$ (см.\ \cite{kp1}) и локальной тривиальности
расслоения $\D^0\to M^{|\tilde\N_{f_*}|}$, $h\mapsto
h(\tilde\N_{f_*})=\tilde\N_{f_*h^{-1}}$ (со слоем $\D^*\cap\D^0$ над
точкой $\tilde\N_{f_*}$, см.\ \cite{Birman}), существуют окрестность
$\UU\subset\FF^1$ оснащенной функции Морса $(f_0,\alpha_0)$ в $\FF^1$
и непрерывное отображение $H\:\UU\to\D^0$, такие что
$H(f_0,\alpha_0)=h_0$ и $(H(f,\alpha))^*(f,\alpha)\in(\FF^*)^1$ для
любой $(f,\alpha)\in\UU$. Получаем непрерывное отображение
$\UU\to(\FF^*)^1/(\D^*\cap\D^0)$,
$(f,\alpha)\mapsto(\D^*\cap\D^0)((H(f,\alpha))^*(f,\alpha))$,
совпадающее с композицией
$\UU\hookrightarrow\FF^1\to\FF^1/\D^0\stackrel{(\overline
i)^{-1}}{\longrightarrow}(\FF^*)^1/(\D^*\cap\D^0)$, откуда следует
непрерывность отображения $(\overline i)^{-1}$.

Пусть $|\N^*|\le\chi(M)+1$. Включение $\D^*\cap\D^0\supseteq(\D^*)^0$
очевидно.
Покажем, что $\D^*\cap\D^0\subseteq(\D^*)^0$. Пусть $\T^*\subset\D^*$
-- подгруппа, аналогичная $\T\subset\D$, см.\ обозначение~\ref
{not:R:G0*}(B). Так как $\D^*\cap\D^0\subseteq\D^*\cap\T=\T^*$ и
количество фиксированных точек $|\N^*|\le\chi(M)+1$, то из (\ref
{eq:EE1}) следует $(\D^*)^0=\T^*\supseteq\D^*\cap\D^0$.
Лемма доказана.
 \end{proof}

\begin{Sta} \label{sta:MM3}
В случае {\rm(\ref {eq:main})} правое действие группы $\T\subset\D$
{\rm(см.\ обозна\-че\-ние~\ref {not:R:G0*})} на $\FF^1$ является
свободным. Имеется $\D^0$-эквивариантный гомеоморфизм
$p_3\:\FF^1\stackrel{\approx}{\longrightarrow}\D^0\times\MM$,
композиция которого с проекцией $\D^0\times\MM\to\MM$ совпадает с
$\Ev$. Здесь группа $\D^0$ действует на $\D^0\times\MM$ справа по
формуле $(h_1,h_2,m)\mapsto(h_2h_1,m)$. В частности, вычисляющее
отображение $\Ev\:\FF^1\to\MM$ является тривиальным
$\D^0$-расслоением, а полный прообраз $\Forg_1^{-1}([f]_{\isot})$ любого класса изотопности $[f]_{\isot}$ в $F^1$ при забывающем
отображении $\Forg_1\:\FF^1\to F^1$ гомеоморфен прямому произведению
 $
 \D^0\times\uups_{[f]_\isot}\approx\D^0\times\MM_{[f]_\isot}=\D^0\times\Ev(\Forg_1^{-1}([f]_{\isot}))
 \subset\D^0\times\MM$.
\end{Sta}

\begin{proof}
{\it Шаг 1.} Докажем свободность действия подгруппы $\T\subset\D$ на
$\FF^1$. Если $h\in\T$ и $h^*(f,\alpha)=(f,\alpha)$, то согласно
\cite[лемма 3.4]{KP3} выполнено $h\in(\stab_{\D^0}f)^0$, а потому $h$ переводит в
себя каждое ребро графа $G_f$ (см.\ обозначение \ref {not:Gf:or}).
Отсюда следует, что $h=\id_M$ ввиду единственности диффеоморфизма
$h\in\D$, сохраняющего оснащенную функцию Морса $(f,\alpha)$ и
переводящего любое ребро графа $G_f$ в себя (см.\ доказательство
утверждения~\ref {sta:MM2}, шаг 3, инъективность).

{\it Шаг 2.} Построим непрерывное отображение $s\:\MM\to\FF^1$, такое что $\Ev\circ s=\id_\MM$ (т.е.\ правое обратное отображения $\Ev$). Согласно утверждению
\ref {sta:MM1} пространство $\MM$ является гладким $3q$-мерным
многообразием и покрыто открытыми подмножествами $\MM_{\succeq[f]_\isot}\subset\MM$.
Фиксируем на $\MM$ клеточное разбиение
(состоящее, вообще говоря, из бесконечного числа клеток), такое что
каждая его замкнутая клетка целиком содержится в одной из областей
$\MM_{\succeq[f]_\isot}$ и характеристическое отображение любой
клетки является гомеоморфизмом на свой образ. Рассмотрим два случая.

{\it Случай 1.} Предположим, что $\chi(M)<p^*+q^*+r^*$.
Пусть $\MM^{(k)}$ -- $k$-мерный остов клеточного разбиения,
$k\le\dim\MM=3q$. Будем строить отображение
$s_{k}\:\MM^{(k)}\to\FF^1$, такое что $\Ev\circ
s_{k}=\id_{\MM^{(k)}}$, индукцией по $k$. При $k=-1$ строить нечего,
так как $\MM^{(-1)}=\emptyset$. Пусть $k\ge0$ и отображение $s_{k-1}$
построено. Рассмотрим любую $k$-мерную клетку
$\sigma=\sigma^k\subset\MM^{(k)}$ разбиения. По построению ее
замыкание целиком содержится в одной из областей
$\MM_{\succeq[f]_\isot}$. Выберем какое-либо поднятие
 $$
\ell_\sigma\:\overline{\sigma}\to S_{\succeq[f]_\isot}\times
U_f^{\infty}
 $$
замкнутой клетки $\overline\sigma$ при накрытии
$p_X|_{\{[f]_\isot\}\times S_{\succeq[f]_\isot}\times
U_f^{\infty}}\circ a_{[f]_\isot}^{-1}\:S_{\succeq[f]_\isot}\times
U_f^{\infty}\to\MM_{\succeq[f]_\isot}$, где
$a_{[f]_\isot}\:\{[f]_\isot\}\times S_{\succeq[f]_\isot}\times
U_f^{\infty}\to S_{\succeq[f]_\isot}\times U_f^{\infty}$ -- проекция.
Тогда
 \begin{equation} \label{eq:*2}
 p_X\circ a_{[f]_\isot}^{-1}\circ \ell_\sigma=\id_{\overline\sigma}.
 \end{equation}
Рассмотрим два $(k-1)$-мерных сфероида в $\FF^1$:
 $$
 S_1:=s_{k-1}|_{\d\sigma^k}\:S^{k-1}\approx\d\sigma^k\to\FF^1,
 \qquad
 S_2:=\i_{f}\circ \ell_\sigma|_{\d\sigma^k}\:S^{k-1}\approx\d\sigma^k\to\FF^1,
 $$
см.\ (\ref {eq:i4f}). Тогда $\Ev\circ S_i=\id_{S^{k-1}}$, $i=1,2$,
так как $\Ev\circ s_{k-1}=\id_{\MM^{(k-1)}}$ в силу индукционного
предположения и
 \begin{equation} \label {eq:*...}
 \Ev\circ\i_f\circ\ell_\sigma=p_X\circ a_{[f]_\isot}^{-1}\circ \ell_\sigma=\id_{\overline\sigma}
 \end{equation}
в силу (\ref {eq:*1}) и (\ref {eq:*2}).
Поэтому (в силу инъективности $\overline{\Ev}$, см.\ утверждение~\ref
{sta:MM2}) для любого $m\in S^{k-1}$ существует диффеоморфизм
$h_m\in\D^0$, такой что $S_1(m)=h_m^*(S_2(m))$. Этот диффеоморфизм
единствен в силу свободности действия группы $\D^0\subset\T$ на
$\FF^1$ (см.\ шаг 1). Получаем однозначное отображение
 $$
 H=H_{\d\sigma^k}\: \d\sigma^k\approx S^{k-1}\to\D^0, \qquad m\mapsto h_m, \quad m\in\d\sigma^k.
 $$
Докажем непрерывность отображения $H$. Так как сфероиды $S_1$ и $S_2$
непрерывны (по индукционному предположению и в силу непрерывности
$\i_f$ и $\ell_\sigma$, см.\ (\ref {eq:i4f})), то они задают
непрерывную зависимость пары оснащенных функций Морса
$S_1(m)=:(f,\alpha)$ и $S_2(m)=:(f_2,\alpha_2)$ от точки $m\in
S^{k-1}$. Если точка $\tilde m\in S^{k-1}$ близка к $m$, то в силу
(\ref {eq:*.}) и (\ref {eq:Ev}) выполнено $S_1(\tilde m)=:(\tilde
f,\tilde\alpha)\stackrel{\Ev}{\mapsto}p_X([f_0]_\isot,h_{f,f_0}^{*0}(\barc(\tilde
fh)),h_{f,f_0}^*[h^*\tilde\alpha])$, $S_2(\tilde m)=:(\tilde
f_2,\tilde\alpha_2)\stackrel{\Ev}{\mapsto}p_X([f_0]_\isot,h_{f_2,f_0}^{*0}(\barc(\tilde
f_2h_2)),h_{f_2,f_0}^*[h_2^*\tilde\alpha_2])$, где $f_0\in F^1$ --
отмеченная функция класса изотопности $[f]_\isot=[f_2]_\isot$,
диффеоморфизмы $h,h_2\in\D^0$ близки к $\id_M$ и $h(\N_f)=\N_{\tilde
f}$, $h_2(\N_{f_2})=\N_{\tilde f_2}$, $h_{f_2,f_0}:=h_mh_{f,f_0}$.
Так как $h_{f,f_0}^*(f,\alpha)=h_{f_2,f_0}^*(f_2,\alpha_2)$,
$\Ev(S_1(\tilde m))=\Ev(S_2(\tilde m))$, то
$(h_{f,f_0}^{*0}(\barc(\tilde
fh)),h_{f,f_0}^*[h^*\tilde\alpha])=(h_{f_2,f_0}^{*0}(\barc(\tilde
f_2h_2)),h_{f_2,f_0}^*[h_2^*\tilde\alpha_2])$, поскольку $p_X$
локально является гомеоморфизмом (см.\ конец \S\ref {sec:MM}). Так
как $\barc(\tilde fh)=\barc(\tilde f_2h_2h_m)$, то согласно критерию
изотопности возмущенных функций Морса (см.\ \cite[утверждение~1.1
и~\S3]{K} или (\ref {eq:incid})) выполнено $\tilde fhh_0=\tilde
f_2h_2h_m$ для некоторого $h_0\in\Diff^0(M,\N_f)$, такого что
автоморфизм $dh_0|_{\b_j(f)}\:T_{\b_j(f)}M\to T_{\b_j(f)}M$ близок к
$\id_{T_{\b_j(f)}M}$ для любой седловой критической точки
$\b_j(f)\in\N_{f,1}$. Отсюда и из равенств $\tilde fhh_0=\tilde
f_2h_2h_m$, $[(hh_0)^*\tilde\alpha]=[(h_2h_m)^*\tilde\alpha_2]$ и
$h_{\tilde m}^*(\tilde f_2,\tilde\alpha_2)=(\tilde f,\tilde\alpha)$
следует (согласно доказательству утверждения \ref {sta:MM1}, шаг 3,
инъективность), что изоморфизм $dh_{\tilde m}|_{\b_j(\tilde
f)}\:T_{\b_j(\tilde f)}M\to T_{\b_j(\tilde f_2)}M$ близок к
изоморфизму $dh_{m}|_{\b_j(f)}\:T_{\b_j(f)}M\to T_{\b_j(f_2)}M$
($1\le j\le q$). Отсюда и из равенства $h_{\tilde m}^*(\tilde
f_2,\tilde \alpha_2)=(\tilde f,\tilde\alpha)$ следует, что в
некоторой окрестности любой седловой критической точки функции
$(f,\alpha)$ выполнено $h_{\tilde m}\to h_m$ при $\tilde m\to m$.
Отсюда и из равенств $h_m^*(df_2^2+\alpha_2^2)=df^2+\alpha^2$,
$h_{\tilde m}^*(d\tilde f_2^2+\tilde\alpha_2^2)=d\tilde
f^2+\tilde\alpha^2$ следует, что $h_{\tilde m}\to h_m$ всюду на $M$
при $\tilde m\to m$. Поэтому диффеоморфизм $h_m$ непрерывно зависит
от $m\in S^{k-1}$, т.е.\ сфероид $H$ непрерывен.

Так как $\D^0=\Diff^0(M,\N)$ и количество фиксированных точек
$|\N|=p^*+q^*+r^*>\chi(M)$, то топологическая группа $\D^0$
стягиваема (см.\ (\ref {eq:EE})), откуда сфероид $H$ непрерывно
продолжается на всю замкнутую клетку $\overline{\sigma^k}$. Пусть
$\tilde H\:\overline{\sigma^k}\to\D^0$, $m\mapsto\tilde h_m$ -- такое
продолжение. Определим отображение $s_k\:\MM^{(k)}\to\FF^1$ формулой
 $$
 s_{k}|_{\overline{\sigma^k}}\:\overline{\sigma^k}\to\FF^1, \qquad m\mapsto \tilde h_m^*(\i_{f}\circ \ell_\sigma(m)).
 $$
Оно однозначно и является продолжением отображения $s_{k-1}$, так как
 $$
 s_{k}|_{\d\sigma^k}\:m\mapsto h_m^*(\i_{f}\circ \ell_\sigma(m))=h_m^*(S_2(m))=S_1(m)=s_{k-1}|_{\d\sigma^k}(m),
 \quad m\in\d\sigma^k.
 $$
При этом $\Ev\circ s_{k}|_{\overline{\sigma^k}}(m)=\Ev(\tilde
h_m^*(\i_{f}\circ \ell_\sigma(m)))=\Ev(\i_{f}\circ \ell_\sigma(m))=m$, $m\in\overline{\sigma^k}$, ввиду (\ref {eq:*...}), откуда
$\Ev\circ s_{k}=\id_{\MM^{(k)}}$. Итак, существование непрерывного отображения
$s$, являющегося правым обратным $\Ev$, доказано в случае $\chi(M)<p^*+q^*+r^*$.

{\it Случай 2.} Предположим теперь, что $p^*+q^*+r^*\le\chi(M)$.
В силу условия (\ref {eq:main}) количество
$|\hat\N_f\setminus\N|=\hat p+\hat q+\hat r-(p^*+q^*+r^*)$
отмеченных, но не фиксированных, критических точек любой функции
$f\in F$ превосходит $\chi(M)-(p^*+q^*+r^*)\ge0$. Поэтому имеется непустое подмножество
$\tilde\N_{f_*}\subseteq\hat\N_{f_*}\setminus\N$, состоящее из
$\chi(M)-(p^*+q^*+r^*)+1>0$
точек. Рассмотрим соответствующие подпространства $\FF^*\subset\FF$ и
$(\FF^*)^1\subset\FF^1$, подгруппы $\D^*\subset\D$ и
$(\D^*)^0\subset\D^0$, и $3q$-мерное многообразие
$\MM^*\approx(\FF^*)^1/(\D^*)^0$, см.\ обозначение \ref {not:fix}.

Так как количество фиксированных точек
$|\N^*|=|\N|+|\tilde\N_{f_*}|=\chi(M)+1>\chi(M)$, то согласно случаю
1 существует непрерывное отображение $s^*\:\MM^*\to(\FF^1)^*$, такое
что $\Ev^*\circ s^*=\id_{\MM^*}$. Так как количество фиксированных
критических точек $|\N^*|=|\N|+|\tilde\N_{f_*}|\le\chi(M)+1$, то по
лемме \ref {lem:KK*} имеется гомеоморфизм
$\overline{i}\:(\FF^*)^1/(\D^*)^0=(\FF^*)^1/(\D^*\cap\D^0)\stackrel{\approx}{\longrightarrow}\FF^1/\D^0$.
Положим
 $$
 s:=i\circ s^*\circ \overline{\Ev^*}\circ(\overline{i})^{-1}\circ{\overline{\Ev}}^{-1}\:\MM\to\FF^1.
 $$
Из определения отображений $\Ev,\Ev^*$ следует, что
 $\Ev|_{\Im\,i}=
 \overline{\Ev}\circ\overline{i}\circ{\overline{\Ev^*}}^{-1}\circ\Ev^*\circ i^{-1}$.
Поэтому $\Ev\circ s=\Ev\circ i\circ s^*\circ
\overline{\Ev^*}\circ(\overline{i})^{-1}\circ{\overline{\Ev}}^{-1}=\id_\MM$.

{\it Шаг 3.} На шаге 2 построено непрерывное отображение
$s\:\MM\to\FF^1$, такое что $\Ev\circ s=\id_\MM$. Определим
непрерывное $\D^0$-эквивариантное отображение
$i_3\:\D^0\times\MM\to\FF^1$ формулой $i_3(h,m):=h^*(s(m))$. Оно
биективно в силу $\Ev\circ s=\id_\MM$, свободности действия $\D^0$ на
$\FF^1$ и биективности $\overline{\Ev}$ (см.\ утверждение \ref
{sta:MM2}). Обратное отображение имеет вид
$p_3=i_3^{-1}\:\FF^1\to\D^0\times\MM$,
$(f,\alpha)\mapsto(\delta(f,\alpha),\Ev(f,\alpha))$, где отображение
$\delta\:\FF^1\to\D^0$ определяется условием
$(\delta(f,\alpha))^*(s\circ\Ev(f,\alpha))=(f,\alpha)$. Его
непрерывность доказывается аналогично доказательству непрерывности
сфероида $H$ (см.\ шаг 2, случай 1). Так как отображения $i_3,p_3$
непрерывны и взаимно обратны, они являются взаимно обратными
гомеоморфизмами. Утверждение \ref {sta:MM3} доказано.
\end{proof}

Утверждение \ref {sta:MM3} доказывает гоме\-о\-мор\-физмы $\FF^1\approx\D^0\times\MM$ и
$\Forg_1^{-1}([f]_{\isot})\approx\D^0\times\MM_{[f]_\isot}$.
С учетом (\ref {eq:EE}) и того, что отображения включения $F^1\hookrightarrow F$, $\FF^1\hookrightarrow\FF$ и забывающие отображения $\FF\to F$,
$\FF^1\to F^1$ и $\Forg_1^{-1}([f]_{\isot})\to[f]_{\isot}$ являются гомотопическими эквивалентностями согласно \cite[теорема 2.5]{kp1}, получаем теорему~\ref
{thm:KP4add}.


\end{document}